\documentclass[a4paper,12pt,leqno]{article}
\usepackage{latexsym}
\usepackage[all]{xy}

\usepackage{amssymb} 
\usepackage{amsmath} 
\usepackage{theorem}

\def\versionnumber{ver.17}

\def\Z{{\mathbb{Z}}}
\def\K{{\mathbb{K}}}
\def\CC{{\mathbb{C}}}
\def\R{{\mathbb{R}}}
\def\A{{\mathcal{A}}}
\def\B{{\mathcal{B}}}

\DeclareMathOperator{\codim}{codim}

\DeclareMathOperator{\Der}{Der}

\DeclareMathOperator{\POexp}{POexp}

\numberwithin{equation}{section}

\newcommand{\owari}{\hfill$\square$}

\theoremstyle{break}
\newtheorem{theorem}{Theorem}[section]
\newtheorem{prop}[theorem]{Proposition}
\newtheorem{cor}[theorem]{Corollary}
\newtheorem{lemma}[theorem]{Lemma}
\newtheorem{define}[theorem]{Definition}

\newtheorem{rem}[theorem]{Remark}

\newtheorem{example}[theorem]{Example}
\newtheorem{problem}[theorem]{Problem}

\title{Roots of the characteristic polynomials of hyperplane arrangements and their 
restrictions and localizations}
\author{Takuro Abe
\footnote
{
Institute of Mathematics for Industry,
Kyushu University,
Fukuoka 819-0395, Japan.
Email:abe@imi.kyushu-u.ac.jp. Tel:+81-92-802-4479.
}
(\versionnumber)
}
\date{\today}

\begin{document}

\maketitle

\begin{abstract}
Terao's factorization theorem shows that if an arrangement is free, then its characteristic polynomial 
factors into the product of linear polynomials over the integer ring. This is not a necessary condition, but 
there are not so many non-free arrangements whose characteristic polynomial factors over the integer ring. On the 
other hand, the localization of a free arrangement is free, and its restriction is in many cases free, thus 
its characteristic polynomial factors. In this paper, we consider how their integer, or real roots behave.
\end{abstract}

\section{Introduction}
Let $\A$ be a central arrangement in 
$V=\K^\ell$ and $H \in \A$. 
A triple $(\A,\A',\A^H)$ is defined by 
$\A':=\A \setminus \{H\}$, and 
$
\A^H:=\{L \cap H \mid L \in \A'\}$. The triple has played a very 
important role in the theory of hyperplane arrangements. Let us show it by example. 
An intersection lattice $L(\A)$ of $\A$ is defined by 
$$
L(\A):=\{\cap_{L \in \B} L \mid \B \subset \A\}
$$
with an order by reverse inclusion. Then $L(\A)$ has a poset structure with 
a M\"obius function $\mu:L(\A) \rightarrow \Z$ defined by 
$\mu(V):=1$, and by 
$$
\mu(X):=-\sum_{X \subsetneq Y \subset V} \mu(Y).
$$
Then a characteristic polynomial $\chi(\A;t)$ is defined by 
$$
\chi(\A;t):=\sum_{X \in L(\A)} \mu(X)t^{\dim X}.
$$
It is known that $(-1)^\ell\chi(\A;-t^{-1})=\mbox{Poin}(V \setminus 
\cup_{L \in \A} L;t)$ when $\K=\CC$. However, this important 
invariant $\chi(\A;t)$ is not easy to compute in general. The most useful 
inductive method is so called the deletion-restriction theorem as follows:
$$
\chi(\A;t)=\chi(\A';t)-\chi(\A^H;t).
$$
For the proof, see Corollary 2.57 in 
\cite{OT} for example. 
We may apply this to compute $\chi(\A;t)$ efficiently. 

Also, the triple plays a key role in algebra of arrangement. Let 
$S:=\mbox{Sym}^*(V^*) \simeq \K[x_1,\ldots,x_\ell]$ be the 
coordinate ring of $V$, and $\Der S =\oplus_{i=1}^\ell S \partial_{x_i}$ the 
module of $\K$-linear $S$-derivations. For a fixed linear form $\alpha_L \in V^*$ of 
$L \in \A$, the logarithmic derivation module $D(\A)$ is defined by 
$$
D(\A):=\{\theta \in \Der S \mid 
\theta(\alpha_H) \in S \alpha_H\ (\forall H \in \A)\}.
$$
$D(\A)$ is a reflexive $S$-module, and we say that $\A$ is free with 
$\exp(\A)=(d_1,\ldots,d_\ell)$ if $D(\A)$ is a free $S$-module with 
a homogeneous basis $\theta_1,\ldots,\theta_\ell$ of degree $\deg \theta_i=d_i$ 
for $i=1,\ldots,\ell$. Here $\deg \theta_i:=\deg \theta_i(\alpha)$ for some linear form 
$\alpha$ with $\theta_i(\alpha) \neq 0$. Freeness has been one of the most important 
property of arrangements, and has been intensively studied by a lot of mathematicians. 
In particular, due to Terao's factorization theorem in \cite{T2} (see Theorem \ref{Teraofactorization}), when $\A$ is free with $\exp(\A)=(d_1,\ldots,d_\ell)$, 
it holds that 
$$
\chi(\A;t)=\prod_{i=1}^\ell (t-d_i).
$$
Hence freeness determines $\chi(\A;t)$, and its roots have algebraic meaning. 
Let us call $\A$ \textbf{integer rooted} if $\chi(\A;t)$ has $\ell$-integer roots, and 
\textbf{real rooted} if $\chi(\A;t)$ has $\ell$-real roots.

It is known that the \textbf{localization} $\A_X:=\{H \in \A\mid X \subset H\}$ of $\A$ at $X \in 
L(\A)$ is free if $\A$ is free. It was conjectured by Orlik that the restriction $\A^H$ of a free arrangement 
$\A$ is free, whose counter example was found by Edelman and Reiner in \cite{ER}, but still almost all cases 
the restriction of a free arrangement is free. By Terao's factorization, this says that 
$\A_X$ (and in many cases $\A^H$) is integer rooted if $\A$ is free. So it is natural to ask the following 
question.

\begin{problem}
Describe the roots of $\chi(\A_X;t)$ and $\chi(\A^H;t)$ when they are both integer, or real rooted and $\A$ is 
integer, or real rooted.
\label{mainproblem}
\end{problem}

One of the known result is the following:

\begin{theorem}[Theorem 1.1, Corollary 1.2, \cite{A}]
Let $\A$ be an arrangement in $\K^3$ and 
$H \in \A,\ \A':=\A 
\setminus \{H\}$. 

(1)\,\, 
If $\chi(\A;t)=(t-1)(t-a)(t-b)$ with 
$1 \le a \le  b$, then $a+1<|\A^H|<b+1
$ cannot occur for any $H \in \A$. Moreover,

(2)\,\,
assume that $\A$ is free, and 
$\chi(\A;t)=(t-1)(t-a)(t-b),\ a
\le b$. Then $|\A^H| \le a+1$ or coincides with $b+1$. In this 
case $\exp(\A^H)=(1,c)$ with either $c=b$ or
$c \le a$.

(3)\,\,Assume that $\A'$ 
is free and 
$\chi(\A';t)=(t-1)(t-a)(t-b),\ a
\le b$. Then $|\A^H| \ge b+1$ or coincides with $a+1$. 
In this 
case $\exp(\A^H)=(1,c)$ with either $c=a$ or
$c \ge b$.
\label{A}
\end{theorem}

The main result in this article is to prove a higher 
dimensional version of Theorem \ref{A} for both restrictions and 
localizations, see Theorems \ref{main100}
and \ref{practical2}. They give us a lot of avoidances of the integer roots
of the restrictions and localizations of hyperplane arrangements which are also 
integer rooted. Though restrictions and localizations are very different, 
the statements of 
these theorems are the same, strangely. 

Also, we investigate the behavior of integer roots for free arrangements. By Terao's factorization
theorem (Theorem \ref{Teraofactorization}), we know that free arrangements are integer rooted. 
If the exponents of the restriction are contained in that of the original one, then the deletion is also free by 
Terao's addition-deletion theorem (Theorem \ref{addition}). However, there are examples of free arrangements 
$\A$ such that $\A^H$ is also free but $\exp(\A^H) \not \subset \exp(\A)$ (e.g., see \cite{ER}). So let us investigate how 
their integer roots are related. 

The organization of this article is as follows. In \S2 we introduce several results used in this atricle. 
In \S3 we study the behavior of integer roots of the restrictions of integer rooted arrangements.
In \S4 we study the behavior of integer roots of the localizations of integer rooted arrangements. 
In \S5 we investigate the behavior of integer roots of $\A$, $\A^H$ and $\A_X$ when all of them are free.
\medskip


\noindent
\textbf{Acknowledgements}. 
The author is grateful to Masahiko Yoshinaga for the advice of the proof 
of Theorem \ref{main100}. 
The author is 
partially supported by 
JSPS KAKENSHI grant number JP20K20880.

\section{Preliminaries}
In this section let us summarize several definitions and results used in this 
article. Let $\A$ be an arrangement of hyperplanes in $V=\K^\ell$. Without any specification, 
we assume that $\A$ is 
irreducible, i.e., there are no decompositions $V=V_1 \oplus V_2$ such that 
$\dim_\K V_i >0$ for $i=1,2$ and there are no $\A_i \subset \A\ 
(i=1,2)$ such that $\A_i$ is a non-empty arrangement in $V_i$ for $i=1,2$ and 
$\A=\A_1 \cup \A_2$.  
Assume that 
every hyperplane $H \in \A$ is defined by a linear form $\alpha_H=0$. Let 
$S=\K[x_1,\ldots,x_\ell]$ be a coordinate ring of $V^*$ and $\Der S:=
\oplus_{i=1}^\ell S \partial_{x_i}$. Then the \textbf{logarithmic vector field} $D(\A)$ of $\A$ is defined as 
$$
D(\A):=\{\theta \in \Der S \mid \theta(\alpha_H) \in S \alpha_H\ (\forall H \in \A)\}.
$$
$D(\A)$ is a reflexive $S$-graded module, and not free in general. We say that $\A$ is \textbf{free} with $\exp(\A)=(d_1,\ldots,d_\ell)$ if 
$D(\A)$ is a free $S$-module of rank $\ell$ with homogeneous basis $\theta_1,\ldots,\theta_\ell$ such that 
$\deg \theta_i=d_i$. If $\cap_{H 
\in \A} H =\{0\}$, the lowest degree element $\theta_1$ can be chosen as the Euler 
derivation $\theta_E=\sum_{i=1}^\ell x_i \partial_{x_i}$ which is always contained in $D(\A)$. Also, for $H \in \A$, define 
$D_H(\A):=\{\theta \in D(\A) \mid \theta(\alpha_H)=0\}$. To investigate their relation, the most fundamental results are as follows:

\begin{prop}[e.g., Theorem 4.37, \cite{OT}]
$\A_X$ is free if $\A$ is free for all $X \in L(\A)$.
\label{locfree}
\end{prop}

\begin{theorem}[e.g., Proposition 4.45, \cite{OT}]
Let $H \in \A$ and 
$\A':=\A \setminus 
\{H\}$. Then 
there is the Euler exact sequence
$$
0 \rightarrow D(\A') \stackrel{\cdot \alpha_H}{\rightarrow} D(\A) 
\stackrel{\rho^H}{\rightarrow} D(\A^H).
$$
Here $\rho^H(\theta)(\overline{f}):=\overline{\theta(f)}$, where 
$f \in S,\ \theta \in D(\A)$ and $\overline{f}$ is the image of $f \in S$ in 
$S/\alpha_H S$. $\rho^H$ 
is called the \textbf{Euler restriction map}.
\label{Euler}
\end{theorem}

\begin{lemma}[Lemma 4.4, \cite{AT}]
$\rho^H$ is generically surjective.
\label{generic}
\end{lemma}

To state the advantage of free arrangements, let us introduce combinatorics and topology of arrangements. 
Let 
$$
L(\A):=\{ \cap_{H \in \B} H \mid \B  \subset \A\}
$$ 
be the 
\textbf{intersection lattice} of $\A$ with a partial order induced from the 
reverse inclusion. The \textbf{M\"{o}bius function} $\mu : L(\A) 
\rightarrow \Z$ is defiend by $\mu(V)=1$, and by 
$\mu(X):=-\sum_{X  \subsetneq Y \subset V} \mu(Y)$ for $L(\A) \ni Y \subsetneq V$. 
Define 
$$
\chi(\A;t):=-\sum_{X \in L(\A)} \mu(X) t^{\dim X},\ 
\pi(\A;t):=\sum_{X \in L(\A)} \mu(X) (-t)^{\codim X}.
$$
Let $\chi(\A;t)=\sum_{i=0}^\ell b_i(\A) t^{\ell-i} (-1)^i$. 
When $\A \neq \emptyset$, it is known that $\chi(\A;t)$ is divisible by 
$t-1$. Define $\chi_0(\A;t):=\chi(\A;t)/(t-1)=
\sum_{i=0}^{\ell-1} b_i^0(\A)t^{\ell-1-i}(-1)^i$. 
 It is known that 
$b_i(\A)$ is the $i$-th Betti number of $V \setminus \cup_{H \in \A} H$ when 
$\K=\CC$. Then we may relate the exponents of free arrangements and the combinatorics and 
topology as follows:

\begin{theorem}[Terao's factorization, \cite{T2}]
Assume that $\A$ is free with $\exp(\A)=(d_1,\ldots,d_\ell)$. Then 
$\chi(\A;t)=\prod_{i=1}^\ell(t-d_i)$.
\label{Teraofactorization}
\end{theorem}

The following inequality is 
the key 
of this article.

\begin{theorem}[$b_2$-inequality and the division theorem, \cite{A2}]
It holds that 
$$
b_2^0(\A) \ge b_2(\A^H)+(|\A^H|-1)(|\A|-|\A^H|-1),
$$
which is equivalent to
$$
b_2(\A) \ge b_2(\A^H)+|\A^H|(|\A|-|\A^H|).
$$
The equality (called the \textbf{$b_2$-equality for $(\A,H))$} holds only if $\A_X:=\{H \in \A \mid H \supset X\}$ 
is free for all $X \in L(\A^H)$ with $\codim_V X=3$. The $b_2$-equality holds true if 
$\chi(\A^H;t) \mid \chi(\A;t)$. Moreover, $\A$ is free if the $b_2$-equality holds, 
and $\A^H$ is free for some $H \in \A$.
\label{division}
\end{theorem}



%

In the application in the next section, the following properties of freeness play key roles.


\begin{prop}[e.g., Proposition 4.41, \cite{OT}]
Let $H \in \A$. Then there is a homogeneous polynomial $B$ of degree $|\A|-1-|\A^H|$ such that, 
for any $\theta \in D(\A \setminus \{H\})$, it holds that 
$$
\theta(\alpha_H) \in (\alpha_H, B).
$$
\label{B}
\end{prop}

\begin{theorem}[\cite{T1}, Terao's addition-deletion theorem]
Let $H \in \A$, $\A':=\A \setminus \{H\}$ and let $\A'':=\A^H$. Then two of the following 
imply the third:
\begin{itemize}
\item[(1)]
$\A$ is free with $\exp(\A)=(d_1,\ldots,d_{\ell-1},d_\ell)$.
\item[(2)]
$\A'$ is free with $\exp(\A')=(d_1,\ldots,d_{\ell-1},d_\ell-1)$.
\item[(3)]
$\A''$ is free with $\exp(\A'')=(d_1,\ldots,d_{\ell-1})$.
\end{itemize}
Moreover, all the three hold true if $\A$ and $\A'$ are both free.
\label{addition}
\end{theorem}

\begin{rem}
As we can see in Theorem \ref{addition}, it could occur that 
$\A$ and $\A^H$ are free, but $\A'$ is not. For example, the arrangement 
introduced in \cite{ER}, called the Edelman-Reiner arrangement, is free with exponents 
$(1,5,5,5,5)$, and there is $H \in \A$ such that $\A^H$ is free and $\exp(\A^H)=(1,3,3,5)$. 
So in this case $\A \setminus \{H\}$ is not free by 
Theorem \ref{addition}
\label{ERrest}
\end{rem}

The surjecrtivity of the Euler restriction map was studied in \cite{A9}, 
one of which is the following.

\begin{theorem}[Free surjection theorem, Theorem 1.13 in \cite{A9}]
Let $H\in \A$ and assume that $\A'
:=
\A \setminus \{H\}$ is free. Then the Euler restriction map 
$\rho^H:D(\A) \rightarrow D(\A^H)$ is 
surjective.
\label{FST}
\end{theorem}

We say that $\A$ is \textbf{SPOG} with exponents $\POexp(\A)=(1,d_2,\ldots,d_\ell)$ and 
\textbf{level} $d$ if there is a minimal free resolution of $D(\A)$ of the following form:
$$
0 
\rightarrow 
S[-d-1] \stackrel{(f_1,\ldots,f_\ell,\alpha)}{\rightarrow} 
\bigoplus_{i=1}^\ell S[-d_i] \oplus S[-d]
\rightarrow D(\A) \rightarrow 0.
$$
Here $\alpha\neq 0$, the degree $d$-derivation of the SPOG arrangement is called the \textbf{level element}, and 
the above set of generators for the SPOG arrangement is an \textbf{SPOG generator}. For details, see \cite{A5}. SPOG arrangements appear naturally as in the following.

\begin{theorem}[Theorem 1.4, \cite{A5}]
Let $\A$ be free and $H \in \A$. Then $\A \setminus \{H\}$ is either free or SPOG.
\label{SPOG}
\end{theorem}

\begin{theorem}[Theorem 5.5, \cite{A5}]
Let $\A'$ be free in $\K^3$ and $H \not \in \A'$. Then $\A:=\A'\cup  \{H\}$ is free or SPOG.  
The former occurs if and only if $|\A^H|-1 \in \exp(\A')$. If the latter occurs and 
$\exp(\A')=(1,d_2,d_3)$, then $\POexp(\A)=(1,d_2+1,d_3+1)$ 
and level 
$|\A^H|-1$.
\label{SPOG2}
\end{theorem}

\section{Integer roots of restrictions}

In this section let us consider the roots of $\chi(\A;t)$. For that, let us introduce the 
following invariants.

\begin{define}
For non-negative integers $a_1,\ldots,a_\ell$, let 
\begin{eqnarray*}
A(a_1,\ldots,a_\ell):&=&\sum_{1 \le i < j \le \ell} |a_i-a_j|^2,\\
B(a_1,\ldots,a_\ell):&=&\sum_{1 \le i < j \le \ell} a_i a_j.
\end{eqnarray*}
\label{def}
\end{define}

What we are interested in 
is the Orlik's conjecture whose counter example was already found in \cite{ER}.
It asserts that the restriction of a free arrangement is also free. It is not true, but 
only few counter examples have been found. 
Hence by Terao's factorization, it is natural to consider when the characteristic polynomial of the 
restricted arrangement factors over 
$\Z$ when $\A$ is not free. Also, it is natural to ask, which integer could be a root of 
characteristic polynomials. We give some answers to it. 

\begin{theorem}
Let $\A$ be a real rooted $\ell$-arrangement with 
$\chi_0(\A;t)=\prod_{i=2}^\ell(t-d_i),\ d_i \in \Z$. For $H \in \A$, assume that $\A^H$ is 
real rooted with 
$\chi_0(\A^H;t)=\prod_{i=2}^{\ell-1} (t-e_i)$. Then 
$$
A(d_2,\ldots,d_\ell) \le A(e_2,\ldots,e_{\ell-1},|\A|-|\A^H|).
$$
\label{main100}
\end{theorem}


\begin{rem}
In terms of Definition \ref{def}, we may describe Theorem \ref{A} as 
$$
A(a,b) \le A(|\A^H|-1,|\A|-|\A^H|).
$$
\end{rem}

\noindent
\textbf{Proof of Theorem \ref{main100}}. 
First note that 
$$
B(1,d_2,\ldots,d_\ell)=b_2(\A)
$$
and 
$$
B(d_2,\ldots,d_\ell)=b_2^0(\A).
$$
Assume that $A(a_2,\ldots,a_\ell) >A(b_2,\ldots,b_\ell)$ for positive integers 
$a_i,b_j$ such that 
$$
a_2+\cdots+a_\ell=b_2+\cdots+b_\ell.
$$
Then  
by the computation 
\begin{eqnarray*}
2B(a_2,\ldots,a_\ell)-2B(b_2,\ldots,b_\ell)&=&
-A(a_2,\ldots,a_\ell)+A(b_2,\ldots,b_\ell)\\
&\ &+
c(a_2+\cdots+a_\ell)^2-c(b_2+\cdots+b_\ell)^2\\&<&0,
\end{eqnarray*}
here $c$ is non-zero constant when $\ell\ge 4$. Thus $B(a_2,\ldots,a_\ell)<B(b_2,\ldots,b_\ell)$. 
Now we show the statement by induction on $\ell$. When $\ell=3$, Theorem \ref{A} completes the proof. Assume that $\ell >3$. 
Assume that 
$$
A(d_2,\ldots,d_\ell) > A(e_2,\ldots,e_{\ell-1},|\A|-|\A^H|).
$$
Then 
$$
b_2^0(\A)=B(d_2,\ldots,d_\ell) < B(e_2,\ldots,e_{\ell-1},|\A|-|\A^H|)=b_2(\A^H)+(|\A^H|-1)(|\A|-|\A^H|).
$$
This contradicts Theorem \ref{division}. \owari

\medskip

In practice, the following is convenient to apply Theorem \ref{main100}.

\begin{cor}
Let $\A$ be a real-rooted arrangement with $\chi_0(\A;t)=\prod_{i=2}^\ell(t-d_i)$. 
Assume that $d_{\ell-1}+1 \le d_\ell-1$. Then for 
$e_i:=d_i\ (i \le \ell-2),\ e_{\ell-1}:=d_{\ell-1}+1,\ e_\ell:=d_\ell-1$, there are no $H \in \A$ 
such that 
$$
\chi_0(\A^H;t)=\prod_{i=2}^\ell (t-e_i)/(t-e_j)
$$
for all $j=2,\ldots,\ell$.
\label{practical}
\end{cor}

\noindent
\textbf{Proof}. Immediate from the proof of Theorem \ref{main100}.\owari
\medskip

\begin{example}
Assume that $\chi_0(\A;t)=(t-2)(t-4)(t-7)$. Then 
Theorem \ref{main100}, or Corollary \ref{practical} says that $\chi_0(\A^H;t)=(t-e_1)(t-e_2)$ cannot occur if 
$(e_1,e_2)$ is one of 
$$
(2,5),(2,6),(3,3),(3,4),(3,5),(3,6),(3,7),(4,4),(4,5),(4,6),(5,5),(5,6).
$$
\label{ex1}
\end{example}

\section{Integer roots of the localizations}
First let us prove the following $b_2$-inequality for the 
localizations, which was essentially proved implicitly, 
and 
used in the proof of Proposition 4.2, \cite{A3}. 

\begin{theorem}[$b_2$-inequality for localizations]
Let $\A$ be an essential arrangement and let $X \in L_{\ell-1}(\A)$.  
Let $\A_X^e$ be the essenatialization of $\A_X$, i.e., $\A_X =\A_X^e \oplus \emptyset_1$. 
Then 
$$
b_2(\A) \ge b_2(\A_X^e)+|\A_X|(|\A|-|\A^e_X|).
$$
Moreover, the equality holds if and only if $\A$ is strictly lineary fibered with respect to $\A_X$.
\label{SS}
\end{theorem}

\noindent
\textbf{Proof}.
$\A$ can be obtained by adding hyperplanes to $\A_X$ which do not contain $X$. Since every hyperplane 
$H \in \A \setminus \A_X$ intersects with $L \in \A_X$ in such a way that 
$H \cap L \neq H \cap L'$ for all distinct $L,L' \in \A_X$, it holds that 
$$
b_2(\A) \ge b_2(\A_X^e)+|\A_X|(|\A|-|\A^e_X|).
$$
The equality case follows immediately by definition.\owari
\medskip

In the previous section, we considered the roots of the restriction of a free arrangement. 
In this section, we consider those of the localization. It is known that $\A_X$ is free if $\A$ is free by Proposition \ref{locfree}, but the relation between $\exp(\A)$ and $\exp(\A_X)$ were not known. This is contrary to the restriction case, i.e., $\A^H$ may not be free even if $\A$ is free, but the relation between $\exp(\A)$ and $\exp(\A^X)$ are clearer than the localization cases by 
Theorem \ref{addition}. We prove the following.

\begin{theorem}
Let $\A$ be a real rooted $\ell$-arrangement with 
$\chi_0(\A;t)=\prod_{i=2}^\ell(t-d_i),\ d_i \in \R$. For $X \in L_{\ell-1}(\A)$, assume that $\A_X$ is 
real rooted with 
$\chi_0(\A_X^e;t)=\prod_{i=2}^{\ell-1} (t-e_i)$. Then 
$$
A(d_2,\ldots,d_\ell) \le A(e_2,\ldots,e_{\ell-1},|\A|-|\A_X|).
$$
\label{main100100}
\end{theorem}

\noindent
\textbf{Proof}. Apply the same argument as in Theorem \ref{main100}.\owari
\medskip

\begin{theorem}
Let $\A$ be an essential arrangement with $\chi_0(\A;t)=\prod_{i=2}^\ell(t-d_i),\ 
d_i \in \R$.
Assume that $d_{\ell-1}+1 \le d_\ell-1$. Then for 
$e_i:=d_i\ (i \le \ell-2),\ d_{\ell-1}:=d_{\ell-1}+1,\ e_\ell:=d_\ell-1$, there are no $X \in L_{\ell-1}(\A)$ 
such that 
$$
\chi_0(\A_X^e;t)=\prod_{i=2}^\ell (t-e_i)/(t-e_j)
$$
for all $j=2,\ldots,\ell$.
\label{practical2}
\end{theorem}

\noindent
\textbf{Proof}. 
Assume that we have such $X$. Compute 
\begin{eqnarray*}
b_2(\A) &\ge& b_2(\A_X^e)+|\A_X|(|\A|-|\A_X|)\\
&=&B(e_1,\ldots,\hat{e}_j,\ldots,e_\ell)+e_j(
e_1+\cdots+\hat{e}_j+\cdots+e_\ell)\\
&=&B(e_1,\ldots,e_\ell)\\
&>& B(d_1,\ldots,d_\ell)=b_2(\A),
\end{eqnarray*}
a contradiction.\owari
\medskip

\begin{example}
Assume that $\chi_0(\A;t)=(t-2)(t-4)(t-7)$. Then 
Theorem \ref{practical2} says that $\chi_0(\A_X;t)=(t-e_1)(t-e_2)$ cannot occur if 
$(e_1,e_2)$ is one of 
$$
(2,5),(2,6),(3,3),(3,4),(3,5),(3,6),(3,7),(4,4),(4,5),(4,6),(5,5),(5,6)
$$
for $X \in L_3(\A)$.
\label{ex2}
\end{example}

When $\ell=3$, we know that $\chi(\A;t)$ itself restricts the multiplicity of 
intersection points, which can be regarded as a localization version of 
Theorem \ref{A}. 

\begin{cor}
Let $\ell=3$ and $\chi_0(\A;t)=(t-d_1)(t-d_2)$ with $d_1 \le d_2$. Then 
there are no $ p \in L_2(\A)$ such that $d_1 < \mu(p) < d_2$. 
\label{multiple}
\end{cor}

\noindent
\textbf{Proof}. Clear by Theorem \ref{practical2}.\owari
\medskip

\begin{example}
Let $\A$ have 
$$
\chi_0(\A)=t^2-8t+13.
$$
Then this has two real roots $4 \pm \sqrt{3}$. So Corollary \ref{multiple} 
implies that $\A$ cannot have $ p\in L_2(\A)$ such that $\mu(p)=3,4,5$. 
For example, this $\chi_0(\A;t)$ can be attained by the arrangement defined by 
$$
xz(x^2-y^2)(x^2-4y^2)(x^2-9y^2)(y-z)=0.
$$
In this csae $\mu(p)=1$ or $6$ as Corollary \ref{multiple} says.
\end{example}

\section{Restrictions and localizations of free arrangements}

In this section let us focus on the case of free arrangements. 
By 
Theorem 
\ref{Teraofactorization}, $\chi(\A;t)$ has only integer roots if $\A$ is free, and this is 
the most important case of such examples. Let us investigate the relation between 
integer roots and freeness as follows:

\begin{prop}
\begin{itemize}
\item[(1)]
Let $\A$ be a free $\ell$-arrangement with 
$\chi_0(\A;t)=\prod_{i=2}^\ell(t-d_i),\ d_i \in \Z$. For $H \in \A$, assume that $\A^H$ is free with 
$\chi_0(\A^H;t)=\prod_{i=2}^{\ell-1} (t-e_i)$. Let 
$1 \le d_2 \le \cdots \le d_\ell$ and $1 \le e_2 
\le \cdots \le e_{\ell-1}$. Then $\chi(\A^H;t) \mid 
\chi(\A;t)$, or $d_2 \ge e_2$, 
$d_i \ge e_{i-1}$ for all $i \ge 3$.
\item[(2)]
Let $\A':=\A \setminus \{H\}$ be a free $\ell$-arrangement with 
$\chi_0(\A';t)=\prod_{i=2}^\ell(t-d_i),\ d_i \in \Z$. Moreover, assume that $\A^H$ is free with 
$\chi_0(\A^H;t)=\prod_{i=2}^{\ell-1} (t-e_i)$. Let 
$1 \le d_2 \le \cdots \le d_\ell$ and $1 \le e_2 
\le \cdots \le e_{\ell-1}$. Then $d_i \le e_i$, 
for $i=2,\ldots,\ell-1$.
\end{itemize}
\label{freeroot}
\end{prop}

\noindent
\textbf{Proof}.
(1)\,\,If $\chi(\A^H;t) \mid 
\chi(\A;t)$, then Theorem \ref{addition} shows that $\A'$ is free. Thus 
$\chi(\A^H;t) \nmid 
\chi(\A;t)$ implies that $\A'$ is not free.
Let $\theta_E,\theta_2,\ldots,
\theta_\ell$ and $\theta_E,\varphi_2,\ldots,\varphi_{\ell-1}$ be a basis for $D(\A)$ and $D(\A^H)$ 
such that
$\deg \theta_i=d_i$ and $\deg \varphi_i=e_i$. 
Assume that $d_2<e_2$. Then we may assume that, by using the Euler exact sequence 
$$
0 \rightarrow D(\A') \stackrel{\cdot 
\alpha_H}{\rightarrow}
D(\A)
\stackrel{\rho^H}{\rightarrow}D(\A^H),
$$
$\alpha_H \mid \theta_2$. Thus $\theta_E,\theta_2/\alpha_2,\theta_3,
\ldots,\theta_\ell$ form a basis for $D(\A')$. Hence $\A'$ is free, a contradiction. So $d_2 \ge e_2$. 
Assume that $d_i <e_{i-1}$ for some $i \ge 3$. Then $\rho^H(\theta_E),\ldots,\rho^H(\theta_i)$ can be 
expressed as linear combinations of $\rho^H(\theta_E),\varphi_2,\ldots,\varphi_{i-2}$. Since 
$\rho^H$ is generically surjective by Lemma \ref{generic}, it holds that 
$$
i-1\le \mbox{rank}_{S/\alpha_H S} \langle \rho^H(\theta_E),\ldots,\rho^H(\theta_i)\rangle_{S/\alpha_H S}\le i,
$$
which is 
impossible since 
$$
\mbox{rank}_{S/\alpha_H S}\langle\rho^H(\theta_E),\varphi_2,\ldots,\varphi_{i-2}\rangle_{S/\alpha_H S}=i-2
$$
and the inclusion above. 

(2)\,\,
Let $\theta_1,\ldots,\theta_\ell$ be a basis for $D(\A')$ with $\deg 
\theta_i=d_i$. Now recall that $\rho^H$ is surjective by Theorem \ref{FST}. So there are 
$\varphi_1,\ldots,\varphi_{\ell-1}$ in $D(\A)$ such that their images by $\rho^H$ go to the basis for 
$D(\A^H)$, and $\deg    \varphi_i=e_i$. Note that $\varphi_1,\ldots,\varphi_{\ell-1}$ are $S$-independent 
since their images are so. Note that $\varphi_i \in D(\A')$. Assume that $e_k<d_k$ for some $k$. Then 
$$
\langle \varphi_1,\ldots,\varphi_k\rangle_S 
\subset 
\langle \theta_1,\ldots,\theta_{k-1}\rangle_S,
$$
which is impossible. \owari 
\medskip

Let us consider Theorem \ref{A} by using Proposition \ref{freeroot}. 
If $\chi(\A^H)\nmid \chi(\A;t)$, then Proposition \ref{freeroot} says that $\exp(\A^H)=(1,|\A^H|-1)$ with $|\A^H|-1 \le b$.
So Proposition \ref{freeroot} gives another proof of Theorem \ref{A} (2). 

\begin{rem}
In general the statement in Proposition \ref{freeroot} is not true if $\A$ or $\A'$ is not free.
For example, let 
$$
Q(\A')=x(y-z)(x^2-y^2)(x^2-9y^2)z=0
$$
in $\R^3$, which is not free 
but $\chi(\A';t)=
(t-1)(t-3)^2$. 
Let $H:y=0$ and $\A:= \A' \cup \{H\}$, which is free with exponents $(1,2,5)$ and 
$\exp(\A^H)=(1,1)$. So Proposition \ref{freeroot} (2) 
does not hold.
\label{rem1}
\end{rem}

Let us investigate the case in Remark \ref{rem1} for details. 
If $\A$ (resp. $\A'$) is free and $\A'$ (resp. $\A$) is integer rooted, then we have the following avoidance 
for roots.

\begin{theorem}
Let $\A$ be an 
arrangement in $\K^3$, $H \in \A$ and 
$\A':=\A 
\setminus  
\{H\}$. Assume that $\A$ and $\A'$ are both 
integer rooted.
\begin{itemize}
\item[(1)]
Assume that $\A$ is free with $\exp(\A)=(1,a,b),\ 
a\le b$ and $\chi_0(\A';t)=
(t-c)(t-d),\ c \le d$. Then 
$a-1\le c\le d \le b$.
\item[(2)]
Assume that $\A'$ is free with $\exp(\A)=(1,a,b),\ 
a\le b$ and $\chi_0(\A;t)=
(t-c)(t-d),\ c \le d$. Then 
$a\le c\le d \le b+1$.
\end{itemize}
\label{3roots}
\end{theorem}

\noindent
\textbf{Proof}. 
(1)\,\, By Theorem \ref{addition}, the statement is clear if $\A'$ is free. 
Therefore, we may assume that 
$\A'$ is not free. Assume that $c<a$. Since $a+b=c+d+1$, let 
$c=a-k,\ d=b+k-1$ for some $k \ge 2$. By using the 
deletion-restriction formula in \S1, compure 
$$
|\A^H|=b_2(\A)-b_2(\A')=ab-(a-k)(b+k-1)+1=k^2+(b-a-1)k+a+1.
$$
By Theorem \ref{A} (2), $|\A^H| \le a$ since $\A$ is free and $\A'$ is 
not free. Thus 
$$
k(k+(b-a-1))< k(k+(b-a-1))+1\le 0.
$$
If $b\ge a+1$, then this holds if and only if $k=0$ and $b=a+1$, impossible since $k>1$. 
Assume that $a=b$. Then $k(k-1) \le 0$ 
cannot occur if $k >1$.

(2)\,\,
The same argument as (1). \owari
\medskip

About the relation between roots of $\chi(\A;t)$ and $\chi(\A^H;t)$ 
when both are free, we can state a stronger statement than Proposition \ref{freeroot} as 
follows:

\begin{theorem}
Let $\A \ni H$, and assume that 
both $\A$ 
and 
$\A^H$ are free with $\exp(\A)=(1,d_2,\ldots,d_\ell)_\le$ and 
$\exp(\A^H)=(1,e_2,\ldots,e_{\ell-1})_\le$. Assume that $\A'$ is not free, 
Then 
there is some $k,\ 
3 \le k \le \ell-1$ such that 
$e_j \le d_j$ for all $j \le k-1,\ 
d_k<e_k$, and $e_i=d_{i+1}$ for all $i \ge k$. In the case of (ii), if $\{\theta_i\}_{i=1}^\ell
$ is a basis for $D(\A)$ with $\deg 
\theta_i=  d_i$, and 
$\{\varphi_i\}_{i=1}^{\ell-1}
$ a basis for $D(\A^H)$ with $\deg 
\varphi_i=  e_i$, then 
$\rho^H(\theta_{i+1})=\varphi_i$ for $i \ge 
k$. 
\label{conjans}
\end{theorem}

\noindent
\textbf{Proof}.
Assume that $e_j \le d_j$ for $j \le k$, and $d_k<e_k$.
Let $\theta_E,\theta_2,\ldots,\theta_\ell$ be a basis for $D(\A)$ with 
$\deg \theta_i=d_i$, and 
$\rho^H(\theta_E),\varphi_2,\ldots,\varphi_{\ell-1}$ be a basis for $D(\A^H)$ with 
$\deg \varphi_i=e_i$.
Then for the Euler restriction map 
$\rho:=\rho^H:D(\A) 
\rightarrow D(\A^H)$, the images $\rho(\theta_2),\ldots,\rho(\theta_k)$ can be 
expressed as linear combinations of $\varphi_2,\ldots,\varphi_{k-1}$ after appropriate 
modifications by the Euler derivations. Let 
$$
\rho(\theta_i)=\sum_{j=2}^{k-1} \rho(f_{ij})\varphi_j
$$
for $i=2,\ldots,k$. Since $\rho(\theta_2),\ldots,\rho(\theta_k)$ are not independent, 
there is a linear relation among 
them. Let us look for it.
Let $\rho(f_{ij})=:g_{ij}$. Then by the cofactor matrix theory, we know that 
$$
\sum_{s=2}^{k} h_s\rho(\theta_s)=0,
$$
where 
$$
h_s:=(-1)^s \det (g_{ij})_{2 \le i \le k,\ i \neq s, 2 \le j \le k-1}.
$$
Note that 
$$
\deg h_s=\sum_{i=2}^k d_i - d_s-\sum_{i=2}^{k-1} e_i 
$$
Thus the relation is at degree 
$$
\deg h_s \theta_s =\sum_{i=2}^k d_i - d_s-\sum_{i=2}^{k-1} e_i +d_s=
\sum_{i=2}^k d_i-\sum_{i=2}^{k-1} e_i.
$$
Let $H_s$ be a canonical lift of $h_s$ to $S$. Thus there is 
$\theta \in D(\A')$ such that 
$$
\sum_{i=2}^k H_i \theta_i = \alpha_H \theta.
$$
Here $\deg \theta=\sum_{i=2}^k d_i-\sum_{i=2}^{k-1} e_i-1.
$

Note that 
\begin{eqnarray*}
\deg B&=&
|\A'|-|\A^H|=d-1\\
&=&\sum_{i=2}^\ell d_i-\sum_{i=2}^{\ell-1} e_i-1\\
&=&\deg \theta+\sum_{i=k}^{\ell-1}(d_{i+1}-e_i) \ge \deg \theta.
\end{eqnarray*} 
If $d_{i+1} >e_i$ for some $k \le i \le \ell-1$, then Proposition \ref{B} shows that 
$\theta =\sum_{i=1}^\ell a_i \theta_i \in D(\A)$ for some 
$a_i \in S$. So
$$
\sum_{i=2}^k H_i \theta_i = \alpha_H \theta=\alpha_H \sum_{i=1}^\ell a_i \theta_i,
$$
showing that 
$\alpha_H \mid H_i$ for $2 \le i \le k$, so $
h_i=0$ for all $2 \le i \le k$. Note that the matrix $(g_{ij})_{2 \le i \le k,\ 2 \le j \le k-1}$ is of rank $k-2$ since 
$\rho^H$ is generically surjective. So at least one of its cofactor determinant $h_2,\ldots,h_{k-1}$ 
of size $k-2$ is non-zero, a contradiction.  

Thus $d_{i+1}=e_i$ for all $k \le i \le \ell-1$ and $|\A'|-|\A^H|=\deg \theta$. In this case, 
$\deg \theta=d_k-1+\sum_{i=2}^{k-1}(d_i-e_i)$. Since $d_i 
\ge e_i$ for $2 \le i \le k-1$, $\deg \theta \ge d_k-1$. Assume that 
$\deg \theta=d_k-1$. Then $\chi(\A^H;t)=\chi(\A;t)/(t-d_k)$. Since $\A^H$ is free, the
division theorem shows that $\A'$ is free, a contradiction. So $\deg \theta \ge d_k$. 
The rest parts can be proved by an usual degree argument. 
\owari
\medskip

Next let us consider the behavior of integer roots of free arrangements and their 
localizations. First let us introduce an easy but fundamental lemma.

\begin{lemma}
Let $X \in L(\A)$, 
$i:D(\A) \rightarrow D(\A_X)$ be a canonical inclusion and 
$\eta:D(\A_X) \rightarrow D(\A_X^e)$ be the canonical surjection. Here 
$\A_X^e$ is the essentialization of $\A_X$, i.e., $\A_X=\A_X^e \oplus \emptyset_{\dim X}$. 
Let $F:=\eta\circ i$. Then $F$ is generically surjective.
\label{generic2}
\end{lemma}

\noindent
\textbf{Proof}.
If the point $p \not \in H$ for any $H \in \A$, then $D(\A)_p \simeq 
D(\A_X)_p \simeq S_p^{\ell}$, so clear. Assume that $p \in X$. 
Let $p \in X \setminus \cup_{H \in \A \setminus \A_X}H$, i.e., 
$p$ is a generic point in $X$. We may assume that $X=\{x_1=\cdots=x_k=0\}$. 
Then by the locality of $D(\A)$ (e.g., see \cite{OT}, Example 4,123), it holds that 
$D(\A)_p=D(\A_X)_p$, which completes the proof.\owari
\medskip

The following is a 
localization version of Theorem \ref{conjans}.

\begin{theorem}
Let $\A$ be free with $\exp(\A)=(1,d_2,\ldots,d_\ell)_\le$. Let $X \in L_{\ell-1}(\A)$ and assume that 
$\exp(\A_X)=(0,1,e_2,\ldots,e_{\ell-1})_\le$. Then $d_i \ge e_{i-1}$ for $i \ge 2$. Moreover, if 
$d_i \le e_i$ for $i <k$, and $e_k>d_k$ 
with some $2 \le k \le \ell$, then $e_i=d_{i+1}$ for $i \ge k$.
\label{local}
\end{theorem}

\noindent
\textbf{Proof}. 
$d_i \ge e_{i-1}$ is clear by the same argument as in the previous section. 
Assume that 
$d_i \le e_i$ for $i <k$, and $e_k>d_k$. Let us choose a basis $x_1,\ldots,x_\ell$ in such a way 
that $X=\{x_1=\cdots=x_{\ell-1}=0\}$. 
Let $i:D(\A) \rightarrow D(\A_X)$ be a canonical inclusion and 
$\eta:D(\A_X) \rightarrow D(\A_X^e)$ be the canonical surjection. Here 
$\A_X^e$ is the essentialization of $\A_X$, i.e., $\A_X=\A_X^e \oplus \emptyset_1$, so 
$D(\A_X^e)=D(\A_X)/S\partial_{x_\ell}$. 
Let $F:=\eta\circ i$. Note that $F$ is generically surjective. Now let $\theta_E,\theta_2,\ldots,\theta_\ell$ be 
a homogeneous basis for $D(\A)$ with $\deg \theta_i=d_i$, and 
let $\partial_{x_\ell},F(\theta_E),\varphi_2,\ldots,\varphi_{\ell-1}$ be 
a homogeneous basis for $D(\A_X^e)$ with $\deg \varphi_i=e_i$. Since 
$d_k<e_k$, $F(\theta_E),F(\theta_2),\ldots,F(\theta_k)$ are $S$-dependent since they are expressed as linear 
combinations of $F(\theta_E),\varphi_2,\ldots,\varphi_{k-1}$. Here we do not distinguish the image of $\varphi_i$ by the map 
$\eta$. So there is a minimal degree relation 
$$
\sum_{i=1}^k F(f_i \theta_i)=0
$$
in $D(\A_X^e)$. Let $d$ be the degree of this relation. We can lift up this into that in $D(\A_X)$ into the form 
$$
\sum_{i=1}^\ell f_i \theta_i=f\partial_{x_\ell}.
$$
Note that the left hand side belongs to $D(\A)$, so is $f\partial_{x_\ell}$. 
Let $Q':=Q(\A)/Q(\A_X)$, i.e., $Q'$ is the product of linear forms $\alpha_H$ such that 
$\alpha_H$ needs $x_\ell$ when expressed. So $Q'\partial_{x_\ell} \in D(\A_X)_{|\A|-|\A_X|}$, and it is clear that 
every $g \partial_{x_\ell} \in D(\A)$ is of the form $g'Q'\partial_{x_\ell}$ for some $g' \in S$. Hence 
$$
d \ge |\A|-|\A_X|.
$$
Now by the cofactor matrix theory, we know that 
$$
d \le \sum_{i=2}^k d_i-\sum_{i=2}^{k-1} e_i.
$$
Let us compute 
\begin{eqnarray*}
d &\ge& |\A|-|\A_X|\\
&=& \sum_{i=1}^\ell d_i-\sum_{i=1}^{\ell-1}e_i\\
&=&d+\sum_{i=k}^{\ell-1} (d_{i+1}-e_i) \ge d
\end{eqnarray*}
since $d_{i+1} \ge e_i$ which we have proved just above.  Hence $d_{i+1}=e_i$ for $i \ge k$, and 
$F(\theta_{i+1})=\varphi_i$ follows by the argument on the grading and the generic 
surjectivity of $F$. \owari
\medskip

Also, on the basis for $D(\A)$ and $D(\A_X)$, we have the following.

\begin{cor}
Let $\A$ be free with $\exp(\A)=(1,d_2,\ldots,d_\ell)_\le$. Assume that 
there is $k,\ 2\le k \le \ell$ such that 
$e_i \le d_i\ (i\le k-1)$ and $d_k<e_{k}$. 
Then by Theorem 
\ref{local}, $\exp(\A_X)=(0,1,e_2,\ldots,e_{k-1},d_{k+1},\ldots,d_{\ell})_\le$.
Then we may choose the basis $\partial_{x_\ell},\theta_E,\varphi_2,\ldots,\varphi_{k-1},
\theta_{k+1},\ldots,\theta_\ell$ for $D(\A_X)$ in such a way that 
$\deg \theta_i=d_i\ 
(i \ge k+1)$ and $\theta_{k+1},\ldots,
\theta_{\ell}$ 
is 
a part of basis for $D(\A)$.
\label{baseext}
\end{cor}

\noindent
\textbf{Proof}. 
Let $\theta_E,\theta_2,\ldots,\theta_\ell$ be a basis for $D(\A)$ with $\deg 
\theta_i
=d_i$. By Theorem \ref{local}, we know that $d_{i+1}=e_{i}$ for $ i \ge k$. 
Since the unique relation in $\mbox{Im}(F)$ is at 
degree $|\A|-|\A_X|$ which is smaller than $d_{k+1}=e_k$, and $F$ is generically surjective, 
the basis elements for $D(\A_X^e)$ whose degree is at least $e_k$ are covered by 
$\theta_{k+1},\ldots,\theta_\ell$. Now the usual degree arugment shows that, for the basis 
$\partial_{x_\ell},\theta_E,\varphi_2,\ldots,
\varphi_{\ell-1}$, it holds that $F(\varphi_i)=F(\theta_{i+1})$ for $ i 
\ge k$. Thus $\theta_{i+1}-\varphi_i =f_i \partial_{x_\ell}$ for some $f_i \in S$. 
Since $\partial_{x_\ell}$ belongs to 
a basis for $D(\A_X)$, we may 
arrange that $\theta_{i+1}=\varphi_i$ for 
$i 
\ge k$.\owari
\medskip

\begin{example}
Contrary to the restriction cases, it is easy to construct examples of Theorem 
\ref{local} and Corollary \ref{baseext}. Let $\ell=3$ and $\K=\R$. Consider the arrangement $\A$ defined by 
$$
Q(\A)=xyz(x^2-y^2)(x^2-4y^2)(x-z)(y-z)=0.
$$
This is free with 
$\exp(\A)=(1,3,5)$. Let $X:=\{
x=y=0\} \in L_2(\A)$. Then $\A_X$ is free with exponents 
$(0,1,5)$, so $3<5$. Thus Corollary \ref{baseext} says that 
the derivations of degree $5$ can be chosen as the same 
one. In fact, the basis for $D(\A_X)$ consists of 
$$
\partial_z,\theta_E,y(x^2-y^2)(x^2-4y^2)\partial_y+yz(x+y)(x-z)(x^2-4y^2)
\partial_z.
$$
Since 
$$
z(x-z)(y-z)\partial_z,\theta_E,y(x^2-y^2)(x^2-4y^2)\partial_y+yz(x+y)(x-z)(x^2-4y^2)
\partial_z
$$
form a basis for $D(\A)$, Corollary \ref{baseext} 
is true of course in this case. 
\label{ex9}
\end{example}

From the above, it seems that the similar theory to Theorem \ref{B} could hold for 
the localizations, e.g., $\theta \in D(\A_X)$ can be chosen as in $D(\A)$ if $\deg 
\theta >|\A|-|\A_X|$. Unfortunately, in general, it does not hold.

\begin{example}
Let 
$$
Q(\A'):=xyz(x^2-y^2)(x^2-4y^2)(x-z),
$$
$H:=\ker (2y-3z),\ 
L:=\ker(y-z)$. Let 
$\A:=\A' \cup\{H\}$ and
$\B:=\A' 
\cup \{L\}$. Let $X:=\{x=y=0\}$. 
As seen in Example \ref{ex9}, both $\A'$ abd $\B$ are  free with $\exp(\B)=(1,3,5)$ and 
$\exp(\A')=(1,2,5)$. By Theorem 
\ref{addition}, we may choose the degree $5$ basis derivations commonly, say $
\theta \in D(\B)_5 \subset 
D(\A')_5$. Since $\exp(\A_X')=\exp(\B_X)=\exp(\A_X)=(0,1,5)$ and $2<3<5$, Theorem \ref{local} says that 
the degree $5$-basis derivation for them is the same as $\theta$. By 
Theorem \ref{SPOG2}, $\A$ is SPOG with $\POexp(\A)=(1,3,6)$ and level $6$. Also, if $\theta_E,\varphi,
\theta$ form a basis for $D(\A')$, then by \cite{A5} it is known that 
$\theta_E,\alpha_H\varphi,\alpha_H \theta$ together 
with the level element of the form $f\varphi+g\theta$ form 
a SPOG generator for $D(\A)$. So $\theta \not \in D(\A)$. \owari
\label{Ex999}
\end{example}

Though we have Example \ref{Ex999}, by Corollary \ref{baseext}, it is natural to ask the following:

\begin{problem}
Investigate the polynomial $B$-theory for the localizations, i.e., something like Theorem \ref{B} for 
localizations.
\end{problem}

\end{document}